\documentclass[12pt]{article}
\usepackage[utf8]{inputenc}
\usepackage[british]{babel}
\usepackage{cmap}
\usepackage{lmodern}
\usepackage[T1]{fontenc}

\usepackage{amsmath,amssymb,amsthm}
\usepackage{booktabs,array}
\usepackage[a4paper,top=25mm,bottom=25mm,left=25mm,right=25mm]{geometry}
\usepackage[dvipsnames,svgnames,table]{xcolor}
\usepackage{microtype}
\usepackage{xcolor}
\usepackage{listings}
\usepackage{hyperref} 
\hypersetup{
  colorlinks   = true,    		
  urlcolor     = blue,    		
  linkcolor    = blue,    		
  citecolor    = ForestGreen      
}

\usepackage{url}
\newtheorem{theorem}{Theorem}
\newtheorem{proposition}[theorem]{Proposition}
\newtheorem{lemma}[theorem]{Lemma}
\newtheorem{corollary}[theorem]{Corollary}
\theoremstyle{definition}

\newtheorem{problem}[theorem]{Problem}
\newtheorem{example}[theorem]{Example}
\theoremstyle{remark}
\newtheorem{remark}[theorem]{Remark}
\newcommand{\PP}{\mathbb{P}}
\newcommand{\EE}{\mathbb{E}}
\newcommand{\supp}{\operatorname{supp}}

\DeclareMathOperator{\Cov}{Cov}
\def\keywords{\vspace{.5em} 
{\noindent \textit{Keywords}: }}

\def\AMS{\vspace{.5em} 
{\noindent \textbf{\emph{MSC} class}: }}

\title{Condorcet's nontransitive dice: \\ A minimal realization of a complete ranking profile}
\author{S\'andor Boz\'oki\textsuperscript{1,2,}\thanks{~Corresponding author.
E-mail: \texttt{bozoki.sandor@sztaki.hun-ren.hu}.}
\and \href{https://sites.google.com/view/laszlocsato}{L\'aszl\'o Csat\'o}\textsuperscript{1,2}
\and J\'ozsef Temesi\textsuperscript{2}}
\date{\today}

\begin{document}

\maketitle

\begin{center}
\footnotesize
\textsuperscript{1}~Laboratory on Engineering and Management Intelligence,\\
Research Group of Operations Research and Decision Systems,\\
Institute for Computer Science and Control (SZTAKI), Hungarian Research Network (HUN-REN);\\
1518 Budapest, P.O. Box 63, Hungary.\\[3pt]
\textsuperscript{2}~Institute of Operations and Decision Sciences, \\
Department of Operations Research and Actuarial Sciences,\\
Corvinus University of Budapest (BCE), Budapest, Hungary.
\end{center}

\begin{abstract}
\noindent
We construct three independent finite-support random variables whose six strict ordering probabilities reproduce the complete ranking-frequency profile in a classical example of Condorcet's paradox. The construction has total support size seven, and we prove that no realization of this profile has smaller total support size. Here dice are allowed to be weighted: different outcomes need not have equal probabilities. We then formulate the realization problem for an arbitrary vector of six ordering probabilities. We distinguish constraints on the pairwise comparison probabilities from constraints on the full ranking distribution and give an example showing that pairwise feasibility is not sufficient. For the fixed interlaced supports $\{1,4,7\}$, $\{2,5,8\}$, and $\{3,6,9\}$, we derive the polynomial realization system. On the branch where the smallest outcome of the second variable has probability zero, we give a compact polynomial characterization for profiles with their first two ordering probabilities positive. This yields an explicit semialgebraic sufficient condition for the general realization problem and a rational construction, both with elementary proofs.
\end{abstract}


\AMS{\dots}

\section{Introduction}

A nontransitive triple of dice consists of three mutually independent random
variables $X,Y,Z$ satisfying
\[
 \PP(X>Y)>\frac12,\qquad \PP(Y>Z)>\frac12,\qquad
 \PP(Z>X)>\frac12.
\]
This phenomenon has been studied since Steinhaus and Trybu\l a
\cite{SteinhausTrybula1959} and Trybu\l a \cite{Trybula1961}.
Other treatments and related constructions include Moon and Moser
\cite{MoonMoser1967}, Gardner \cite{Gardner1970}, Savage \cite{Savage1994},
Boz\'oki \cite{Bozoki2014}, Conrey et al.\ \cite{Conrey2016}, Buhler,
Graham, and Hales \cite{Buhler2018}, and Komisarski \cite{Komisarski2021}.
The terms \emph{nontransitive dice} and \emph{intransitive dice} are both
used in this literature. Related intransitive devices are discussed by
Poddiakov \cite{Poddiakov2018}.
While nontransitivity might appear to be a mere mathematical curiosity at first sight, nontransitive patterns can be observed among top clubs in European football \cite{vanOurs2024, vanOurs2025, vanOurs2026}. 

Throughout this paper, a \emph{die} means a real-valued random variable with
finite support. Its outcomes need not be equiprobable; equivalently, we allow
\emph{weighted dice}. A \emph{face} is a distinct value with strictly positive
probability. Thus, the number of faces is the support size, not the number of
physical faces in an equiprobable implementation with repeated labels.
All dice in a realization are mutually independent. We minimize the quantity
\[
 N(A,B,C)=|\supp(A)|+|\supp(B)|+|\supp(C)|.
\]
In particular, the minimality statements below concern weighted finite-support
realizations; they do not assert minimality for physical dice with equiprobable
faces.

The connection with majority cycles in voting is classical; see Condorcet
\cite{Condorcet1785}, Gehrlein \cite{Gehrlein2006}, and Young \cite{Young1995}, among others.
Pairwise comparison probabilities, however, do not determine the probabilities
of all six rankings. Our focus is the stronger requirement of reproducing a
\emph{complete ranking distribution}. We realize a classical 60-voter profile (Section~\ref{Sec2}) with total support size seven (Section~\ref{Sec3}) and prove that this is optimal (Section~\ref{Sec4}). We then study necessary conditions for general ranking distributions (Section~\ref{Sec5}) and obtain a compact constructive criterion on a particular interlaced support branch (Section~\ref{Sec6}). Finally, Section~\ref{Sec7} concludes.

\section{A realization of the classical Condorcet profile} \label{Sec2}

Consider the following preference profile of 60 voters, a classical example
associated with Condorcet's paradox; see \cite{Condorcet1785, Cook2006, Gehrlein2006}:
\begin{center}
\begin{tabular}{lc}
\toprule
Ranking & Number of voters\\
\midrule
$A\succ B\succ C$ & 23\\
$A\succ C\succ B$ & 0\\
$B\succ A\succ C$ & 2\\
$B\succ C\succ A$ & 17\\
$C\succ A\succ B$ & 10\\
$C\succ B\succ A$ & 8\\
\bottomrule
\end{tabular}
\end{center}
Here $A$ defeats $B$ by $33$ votes to $27$, $B$ defeats $C$ by $42$ votes to $18$,
and $C$ defeats $A$ by $35$ votes to $25$. We seek independent random variables
that realize the six normalized ranking frequencies, not merely these three
pairwise comparison probabilities.

\begin{proposition}\label{prop:construction}
Let $A,B,C$ be independent with distributions
\begin{align*}
 &\PP(A=1)=\frac{333}{828},\quad
 \PP(A=3)=\frac{35}{828},\quad
 \PP(A=5)=\frac{460}{828},\\
 &\PP(B=4)=\frac{99}{100},\quad \PP(B=7)=\frac1{100},\\
 &\PP(C=2)=\frac{23}{33},\quad \PP(C=6)=\frac{10}{33}.
\end{align*}
Then
\begin{equation}\label{eq:profile}
\begin{aligned}
 \PP(A>B>C)&=\frac{23}{60}, & \PP(A>C>B)&=0,\\
 \PP(B>A>C)&=\frac{2}{60}, & \PP(B>C>A)&=\frac{17}{60},\\
 \PP(C>A>B)&=\frac{10}{60}, & \PP(C>B>A)&=\frac{8}{60}.
\end{aligned}
\end{equation}
\end{proposition}
\begin{proof}
There are $3\cdot2\cdot2=12$ elementary outcomes. The following table assigns
each of them to its strict ordering. Within an entry such as $7>3>2$, the
labels refer to the variables in the order specified in the first column.
\begin{center}
\begin{tabular}{p{25mm}p{91mm}}
\toprule
Ordering & Elementary outcomes\\
\midrule
$A>B>C$ & $5>4>2$\\
$A>C>B$ & None\\
$B>A>C$ & $4>3>2$, $7>3>2$, $7>5>2$\\
$B>C>A$ & $4>2>1$, $7>2>1$, $7>6>1$, $7>6>3$, $7>6>5$\\
$C>A>B$ & $6>5>4$\\
$C>B>A$ & $6>4>1$, $6>4>3$\\
\bottomrule
\end{tabular}
\end{center}
By independence, the probability of each elementary outcome is the product of
the corresponding atom probabilities. Summing within each row gives
\eqref{eq:profile}. For example,
\[
 \PP(A>B>C)=\frac{460}{828}\frac{99}{100}\frac{23}{33}
 =\frac{23}{60},\qquad
 \PP(C>A>B)=\frac{10}{33}\frac{460}{828}\frac{99}{100}
 =\frac{10}{60}.
\]
All labels on different dice are distinct, so ties are impossible.
\end{proof}

\section{A parametric route to the construction} \label{Sec3}

Introduce independent variables with interlaced potential supports:
\begin{equation}\label{eq:supports}
\begin{array}{c|ccc}
 &\text{value}&\text{value}&\text{value}\\ \hline
 X&1\ (p)&4\ (q)&7\ (1-p-q)\\
 Y&2\ (r)&5\ (s)&8\ (1-r-s)\\
 Z&3\ (t)&6\ (u)&9\ (1-t-u)
\end{array}
\end{equation}
The quantities in parentheses are the corresponding probabilities. The
parameters satisfy
\begin{equation}\label{eq:simplices}
 p,q\geq0,\ p+q\leq1;\qquad
 r,s\geq0,\ r+s\leq1;\qquad
 t,u\geq0,\ t+u\leq1.
\end{equation}
These supports provide a convenient search model, not an assumption that
every ranking distribution has a realization of this form.

Set $A=X$, $B=Z$, and $C=Y$. One obtains the target profile from the family
\begin{equation}\label{eq:family}
\begin{gathered}
 t=0,\qquad u=\frac{99}{100},\qquad
 r=\frac{23}{33}-s,\\
 p=\frac{111-176s}{12(23-33s)},\qquad
 q=\frac{35}{36(23-33s)},\qquad
 0\leq s\leq\frac{111}{176}.
\end{gathered}
\end{equation}
This family can be verified directly using the ordering polynomials in Section~\ref{Sec6}, with the indicated permutation of variables.
In particular,
\[
 1-p-q=\frac59,\qquad 1-r-s=\frac{10}{33},\qquad
 1-t-u=\frac1{100}.
\]
The interval in \eqref{eq:family} is exactly the admissible parameter range
for this rational family. Indeed, $s\geq0$ and $r\geq0$ force
$s\leq23/33$; the endpoint $23/33$ is excluded by the denominators. On the
remaining interval, $p\geq0$ is equivalent to $s\leq111/176$.
The other requirements in \eqref{eq:simplices} then follow from the displayed
identities and nonnegativity.

At $s=0$ the nonzero labels are
$1,4,7$ for $A$, $6,9$ for $B$, and $2,8$ for $C$.
Replacing these seven labels, in increasing order, by $1,\ldots,7$ gives
Proposition~\ref{prop:construction}. This is a common order-preserving
relabelling of the union of the supports and hence preserves all comparisons.

At the other endpoint, $s=111/176$, one has
$p=0$, $q=4/9$, and $r=35/528$. An order-preserving relabelling gives a
second realization:
\begin{align*}
 &\PP(A=2)=\frac49,\qquad \PP(A=5)=\frac59,\\
 &\PP(B=4)=\frac{99}{100},\qquad \PP(B=7)=\frac1{100},\\
 &\PP(C=1)=\frac{35}{528},\quad
 \PP(C=3)=\frac{333}{528},\quad
 \PP(C=6)=\frac{160}{528}.
\end{align*}
Thus the total support size is seven at both endpoints of this family and
eight at every interior point.

\section{The minimality of total support size} \label{Sec4}

Throughout this section, a realization means independent finite-support
variables satisfying \eqref{eq:profile}. Since the six strict ordering
probabilities sum to one, every such realization is tie-free. In particular,
no positive-probability atom can be shared by two different dice.
The pairwise probabilities are
\begin{equation}\label{eq:profilemarginals}
 \PP(A>B)=\frac{33}{60},\qquad
 \PP(B>C)=\frac{42}{60},\qquad
 \PP(A>C)=\frac{25}{60}.
\end{equation}

\begin{lemma}\label{lem:constant}
Every die in a realization of \eqref{eq:profile} has at least two support points.
\end{lemma}
\begin{proof}
Suppose that one die $O$ is constant, say $O=m$, and let $D,E$ be the other
two. Independence gives
\[
 \PP(D>O>E)=\PP(D>m)\PP(m>E)
 =\PP(D>O)\PP(O>E).
\]
For $O=A$, the required profile would imply
\[
 \frac{2}{60}=\PP(B>A>C)
 =\frac{27}{60}\frac{25}{60}=\frac3{16},
\]
a contradiction. For $O=B$ it would imply
\[
 \frac{23}{60}=\PP(A>B>C)
 =\frac{33}{60}\frac{42}{60}=\frac{77}{200},
\]
also a contradiction. For $O=C$ it would imply
\[
 0=\PP(A>C>B)=\frac{25}{60}\frac{18}{60}=\frac18,
\]
which is impossible as well.
\end{proof}

\begin{lemma}\label{lem:twopoint}
The profile \eqref{eq:profile} cannot be realized by three dice with two
support points each.
\end{lemma}
\begin{proof}
Write the supports as $a_1<a_2$, $b_1<b_2$, and $c_1<c_2$, with every
atom having positive probability. Let $w=\PP(C=c_1)\in(0,1)$. Independence and
$\PP(A>C>B)=0$ give
\[
 0=\sum_{c\in\{c_1,c_2\}}\PP(C=c)\PP(A>c)\PP(c>B).
\]
Thus $\PP(A>c)\PP(c>B)=0$ for each $c\in\{c_1,c_2\}$.
Since $\PP(A>C)>0$, one of these points, say $c$, satisfies
$\PP(A>c)>0$, and therefore $B>c$ almost surely. Here we use tie-freeness
to exclude $B=c$. Similarly, $\PP(C>B)>0$ gives a point $c'$ with
$\PP(c'>B)>0$, so $A<c'$ almost surely.
These points cannot coincide: at $c=c'$ both factors in the preceding
product would be positive. Moreover, $B>c$ almost surely and
$\PP(B<c')>0$ imply $c<c'$. Hence
\begin{equation}\label{eq:extremebounds}
 B>c_1\quad\text{and}\quad A<c_2\qquad\text{almost surely}.
\end{equation}

Consequently, $A>B>C$ is equivalent to $\{C=c_1,A>B\}$: the case $C=c_2$
is excluded by $A<c_2$, and $B>c_1$ holds automatically in the other case.
Independence of $C$ from $(A,B)$ now yields
\[
 \frac{23}{60}=w\frac{33}{60},\qquad w=\frac{23}{33}.
\]
Using $B>c_1$ again,
\[
 \frac{42}{60}=w+(1-w)\PP(B>c_2),\qquad
 \PP(B>c_2)=\frac1{100}.
\]
This probability lies strictly between zero and one. Since $B$ has two
support points, necessarily $b_1<c_2<b_2$, with
\[
 \PP(B=b_1)=\frac{99}{100},\qquad \PP(B=b_2)=\frac1{100}.
\]
Similarly, $A<c_2$ gives
\[
 \frac{25}{60}=w\PP(A>c_1),\qquad
 \PP(A>c_1)=\frac{55}{92}\in(0,1).
\]
Therefore $a_1<c_1<a_2$ and $\PP(A=a_2)=55/92$.

The event $B>A>C$ can only occur when $C=c_1$ and $A=a_2$: $C=c_2$ is
excluded by \eqref{eq:extremebounds}, while $A=a_1<c_1$ is impossible when
$C=c_1$. Thus it is equivalent to
$\{C=c_1,A=a_2,B>a_2\}$, and independence gives
\[
 \frac{2}{60}=\frac{23}{33}\frac{55}{92}\PP(B>a_2)
 =\frac5{12}\PP(B>a_2).
\]
It follows that $\PP(B>a_2)=2/25$. But a two-point variable $B$ with the
atom probabilities above has, at any threshold, an upper-tail probability
belonging to
\[
 \left\{0,\frac1{100},\frac{99}{100},1\right\}.
\]
None equals $2/25$, giving the required contradiction.
\end{proof}

\begin{theorem} \label{thm:minimal}
The minimum total support size among independent finite-support realizations
of \eqref{eq:profile} is seven.
\end{theorem}
\begin{proof}
Proposition~\ref{prop:construction} gives a realization with total support
size $3+2+2=7$. If a realization had total support size at most six,
Lemma~\ref{lem:constant} would force every die to have exactly two support
points. Lemma~\ref{lem:twopoint} excludes that possibility.
\end{proof}

\section{General ranking distributions} \label{Sec5}

\begin{problem}\label{prob:general}
Given $v=(v_1,\ldots,v_6)$ in the probability simplex
\[
 \Delta_5=\left\{v\in\mathbb{R}_{\geq0}^6:\sum_{i=1}^6v_i=1\right\},
\]
when do there exist independent real-valued random variables $X,Y,Z$ such that
\begin{equation}\label{eq:orders}
\begin{aligned}
 \PP(X>Y>Z)&=v_1, & \PP(X>Z>Y)&=v_2,\\
 \PP(Y>X>Z)&=v_3, & \PP(Y>Z>X)&=v_4,\\
 \PP(Z>X>Y)&=v_5, & \PP(Z>Y>X)&=v_6?
\end{aligned}
\end{equation}
A finite-support realization is a realization by weighted dice.
\end{problem}
Since the six probabilities sum to one, any realization in
Problem~\ref{prob:general} is automatically tie-free. Define its pairwise
comparison probabilities by
\begin{equation}\label{eq:marginals}
 \xi=v_1+v_2+v_5,\qquad
 \eta=v_1+v_3+v_4,\qquad
 \zeta=v_4+v_5+v_6.
\end{equation}
They equal $\PP(X>Y)$, $\PP(Y>Z)$, and $\PP(Z>X)$, respectively.
A realization is nontransitive in this orientation precisely when
$\xi,\eta,\zeta>1/2$; this additional restriction is not part of the general
realization problem.

\subsection{Pairwise feasibility}

The following classical characterization concerns pairwise probabilities,
not the full six-component distribution.
\begin{theorem} \label{thm:trybula}
\cite[Theorem~1]{Trybula1961} 
Let $(\xi,\eta,\zeta)\in[0,1]^3$. There exist independent real-valued random
variables $X,Y,Z$ such that
\[
 \PP(X=Y)=\PP(Y=Z)=\PP(Z=X)=0,
\]
and
\[
 \PP(X>Y)=\xi,\quad \PP(Y>Z)=\eta,\quad \PP(Z>X)=\zeta
\]
if and only if
\begin{equation}\label{eq:alpha}
 1-\alpha(1-\xi,1-\eta)\leq\zeta\leq\alpha(\xi,\eta),
\end{equation}
where
\[
 \alpha(a,b)=
 \begin{cases}
 \displaystyle\max\left\{\frac{1-a}{b},\frac{1-b}{a},1-ab\right\},
       &a+b>1,\\[4pt]
 1, &a+b\leq1.
 \end{cases}
\]
The denominators in the first case are positive.
\end{theorem}

\begin{proof}[Source and equivalent formulation]
This is Trybu\l a's characterization \cite[Theorem~1]{Trybula1961}.
An equivalent symmetric formulation is given in
\cite[Lemma~2.3]{VuksanovicHildebrand2021}, following Suck
\cite{Suck2002}. With $\bar\xi=1-\xi$, $\bar\eta=1-\eta$, and
$\bar\zeta=1-\zeta$, that formulation is
\begin{align}
 \min\{\xi+\eta\zeta,\eta+\zeta\xi,\zeta+\xi\eta\}&\leq1,
       \label{eq:trybula-upper}\\
 \min\{\bar\xi+\bar\eta\bar\zeta,
        \bar\eta+\bar\zeta\bar\xi,
        \bar\zeta+\bar\xi\bar\eta\}&\leq1.
       \label{eq:trybula-lower}
\end{align}
For completeness, these inequalities imply exactly \eqref{eq:alpha}.
If $\xi+\eta\leq1$, inequality~\eqref{eq:trybula-upper} holds for every
$\zeta\in[0,1]$, since $\xi+\eta\zeta\leq\xi+\eta\leq1$.
If $\xi+\eta>1$, at least one of its three terms must be at most one,
which is equivalent to
\[
 \zeta\leq\frac{1-\xi}{\eta}\quad\text{or}\quad
 \zeta\leq\frac{1-\eta}{\xi}\quad\text{or}\quad
 \zeta\leq1-\xi\eta.
\]
Their union is $\zeta\leq\alpha(\xi,\eta)$.
Applying the same argument to the complementary triple yields the lower
bound in \eqref{eq:alpha}. Thus \eqref{eq:alpha} is valid on the whole cube,
including its boundary, not only in the nontransitive region.
\end{proof}

\begin{corollary}\label{cor:pairwise}
Any realization of \eqref{eq:orders} satisfies
\begin{equation}\label{eq:pairwisebounds}
 1\leq\xi+\eta+\zeta\leq2,\qquad
 \xi\eta\zeta\leq\frac14,\qquad
 (1-\xi)(1-\eta)(1-\zeta)\leq\frac14.
\end{equation}
\end{corollary}
\begin{proof}
For three distinct real numbers, the number of true comparisons among
$X>Y$, $Y>Z$, and $Z>X$ is either one or two. Taking expectations gives
the sum bounds; these bounds do not require independence.
By \eqref{eq:trybula-upper}, at least one of
$\xi+\eta\zeta$, $\eta+\zeta\xi$, and $\zeta+\xi\eta$ is at most one.
If, for example, $\xi+\eta\zeta\leq1$, then
\[
 \xi\eta\zeta\leq\xi(1-\xi)\leq\frac14.
\]
The other two cases are identical. Applying the same argument to
\eqref{eq:trybula-lower} gives the complementary product bound.
\end{proof}

\begin{remark}\label{rem:trivial}
Some elementary necessary inequalities do not impose any restriction beyond
$v\in\Delta_5$ and \eqref{eq:marginals}. For example,
$v_1\leq\min\{\xi,\eta,1-\zeta\}$ follows because each term in the minimum
is a sum of nonnegative components containing $v_1$.
Likewise,
\[
 \xi+\eta+\zeta=1+v_1+v_4+v_5\in[1,2]
\]
is automatic for any $v\in\Delta_5$. Such identities should be distinguished
from restrictions caused by independence.
\end{remark}

\subsection{Constraints on the full ranking distribution}
The next elementary bound uses independence and is stronger than the
corresponding event-inclusion bounds.
\begin{lemma}\label{lem:chain}
If $U,V,W$ are independent real-valued random variables, then
\[
 \PP(U>V>W)\leq\PP(U>V)\PP(V>W).
\]
\end{lemma}
\begin{proof}
Let $f(v)=\PP(U>v)$ and $g(v)=\PP(W<v)$. Independence gives
\[
 \PP(U>V>W)=\EE[f(V)g(V)],\quad
 \PP(U>V)=\EE[f(V)],\quad \PP(V>W)=\EE[g(V)].
\]
The function $f$ is nonincreasing and $g$ is nondecreasing. If $V'$ is an
independent copy of $V$, then
\[
 2\Cov(f(V),g(V))
 =\EE\big[(f(V)-f(V'))(g(V)-g(V'))\big]\leq0.
\]
The claimed inequality follows.
\end{proof}

\begin{corollary}\label{cor:chainbounds}
Every realization of \eqref{eq:orders} satisfies
\begin{equation}\label{eq:chainbounds}
\begin{aligned}
 v_1&\leq\xi\eta,&
 v_2&\leq(1-\zeta)(1-\eta),\\
 v_3&\leq(1-\xi)(1-\zeta),&
 v_4&\leq\eta\zeta,\\
 v_5&\leq\zeta\xi,&
 v_6&\leq(1-\eta)(1-\xi).
\end{aligned}
\end{equation}
\end{corollary}
\begin{proof}
Apply Lemma~\ref{lem:chain} to each ordering and use tie-freeness to express
reversed pairwise probabilities by complements.
\end{proof}

\begin{example}\label{ex:pairwisenotenough}
Pairwise feasibility does not suffice for full-profile realizability. Consider
\[
 v=\left(\frac12,0,0,0,0,\frac12\right).
\]
Its pairwise probabilities are $\xi=\eta=\zeta=1/2$, a feasible triple
obtained, for instance, from independent identically distributed continuous
variables. Nevertheless, this particular full profile violates
$v_1\leq\xi\eta$, because $1/2>1/4$. Thus no independent random variables,
discrete or otherwise, realize this $v$. The same profile satisfies the
simplex conditions and all event-inclusion bounds mentioned in
Remark~\ref{rem:trivial}.
\end{example}

\section{The interlaced support model} \label{Sec6}
Return to \eqref{eq:supports}. Put $h=1-p-q$ and enumerate the 27 elementary
outcomes. In the same order as \eqref{eq:orders}, their probabilities are
\begin{equation}\label{eq:polynomials}
\begin{aligned}
 F_1&=hst,\\
 F_2&=(1-p)rt+h(r+s)u,\\
 F_3&=q(1-r)t+h(1-r-s)(t+u),\\
 F_4&=p(1-r)t+(p+q)(1-r-s)u,\\
 F_5&=qr(1-t)+h(r+s)(1-t-u),\\
 F_6&=pr+(p+q)s(1-t)+(1-r-s)(1-t-u).
\end{aligned}
\end{equation}
For example, $X>Y>Z$ occurs only at $(X,Y,Z)=(7,5,3)$, yielding $F_1=hst$.
The other expressions are obtained by collecting the corresponding elementary
products. The six events partition the sample space, so
$F_1+\cdots+F_6=1$ whenever the atom probabilities are normalized.
Consequently, a profile has a realization in this fixed support pattern if
and only if
\[
 F_i=v_i\quad(i=1,\ldots,6)
\]
has a solution satisfying \eqref{eq:simplices}.
This is a polynomial feasibility problem over a product of three probability
simplices, each of dimension two.

\subsection{A compact characterization on the branch $r=0$}

A useful rational branch arises by setting $r=0$, so $Y$ is supported on a
subset of $\{5,8\}$. The following result replaces expanded rational
expressions by a small collection of explicitly defined polynomials.

\begin{theorem}\label{thm:branch}
Let $v\in\Delta_5$ satisfy $v_1>0$ and $v_2>0$. Define the polynomials
\begin{equation}\label{eq:compactpolys}
\begin{aligned}
 d&=v_1+v_2+v_5, & L&=v_1+v_2,\\
 K&=L-d(v_1+v_3+v_4),\\
 Q&=(v_3+L)K-v_2L,\\
 P&=v_1K-v_1v_2d-dQ.
\end{aligned}
\end{equation}
Then $v$ has a realization in \eqref{eq:supports} with $r=0$ if and only if
\begin{equation}\label{eq:criterion}
 0<K\leq v_2,\qquad Q\geq0,\qquad P\geq0.
\end{equation}
When these inequalities hold, the parameters on this branch are uniquely
specified by
\begin{equation}\label{eq:compactsolution}
 p=\frac{P}{v_1K},\qquad
 q=\frac{dQ}{v_1K},\qquad
 r=0,\qquad s=\frac{K}{v_2},\qquad
 t=\frac{v_1}{d},\qquad u=\frac{v_2}{d}.
\end{equation}
In particular,
\begin{equation}\label{eq:remainingatoms}
 h=1-p-q=\frac{v_2d}{K},\qquad
 1-r-s=1-\frac{K}{v_2},\qquad
 1-t-u=\frac{v_5}{d}.
\end{equation}
\end{theorem}
\begin{proof}
\emph{Necessity and derivation.}
Suppose first that a realization with $r=0$ exists. Equations
\eqref{eq:polynomials} reduce to
\begin{equation}\label{eq:rzero}
\begin{aligned}
 v_1&=hst,& v_2&=hsu,& v_5&=hs(1-t-u),\\
 v_3&=qt+h(1-s)(t+u),&&&\\
 v_4&=pt+(1-h)(1-s)u.&&&
\end{aligned}
\end{equation}
The first three equations give $d=hs$. Since $v_1,v_2>0$, all of $h,s,t,u$
are positive, and
\[
 t=\frac{v_1}{d},\qquad u=\frac{v_2}{d},\qquad s=\frac{d}{h}.
\]
The $v_3$ equation gives
\begin{equation}\label{eq:qrelation}
 qv_1=dv_3-(h-d)L=d(v_3+L)-hL.
\end{equation}
Substituting $p=1-h-q$, $s=d/h$, and \eqref{eq:qrelation} into the
$v_4$ equation gives
\begin{align*}
 dv_4
 &= (1-h-q)v_1+(1-h)(1-d/h)v_2\\
 &= L-d(v_1+v_3)-\frac{v_2d}{h}.
\end{align*}
Therefore
\[
 K=L-d(v_1+v_3+v_4)=\frac{v_2d}{h}>0.
\]
It follows that $h=v_2d/K$ and $s=K/v_2$. Since $s\leq1$, we have $K\leq v_2$.
Equation~\eqref{eq:qrelation} now yields $q=dQ/(v_1K)$, and
$p=1-h-q=P/(v_1K)$. The nonnegativity of $p,q$ is equivalent to
$P,Q\geq0$, since $v_1,K,d>0$. This proves necessity as well as the formulas
and uniqueness on this branch.

\emph{Sufficiency.}
Conversely, assume \eqref{eq:criterion} and define the parameters by
\eqref{eq:compactsolution}. All denominators are positive. The defining
identity for $P$ gives
\[
 p+q+\frac{v_2d}{K}=1.
\]
The three terms are nonnegative, so they are valid atom probabilities for $X$.
Likewise, $0<s\leq1$ gives valid probabilities for $Y$, and
\[
 t+u+\frac{v_5}{d}=1
\]
gives valid probabilities for $Z$. In particular, \eqref{eq:simplices} holds.
Since $hs=d$, the first three equations in \eqref{eq:rzero} are immediate.
The formula for $q$ gives \eqref{eq:qrelation}, proving the equation for
$v_3$. Reversing the displayed calculation for $dv_4$ and using the definition
of $K$ proves the equation for $v_4$. The remaining equation, $F_6=v_6$,
follows because both vectors sum to one. No symbolic computation is needed
for this argument.
\end{proof}

\begin{remark}\label{rem:semialgebraic}
All expressions in \eqref{eq:compactpolys} are polynomials in
$v_1,\ldots,v_5$, with $v_6=1-\sum_{i=1}^5v_i$. Thus the simplex constraints,
$v_1,v_2>0$, and \eqref{eq:criterion} are an explicit semialgebraic condition
without any unresolved denominator signs. This condition is necessary and
sufficient on the branch $r=0$ under $v_1,v_2>0$, but only sufficient for
Problem~\ref{prob:general} or for the unrestricted interlaced model.
Other solutions may have $r>0$, and profiles with $v_1v_2=0$ require separate
analysis. The original Condorcet profile lies on such a boundary in its
original ordering; its realization is supplied separately in
Section~\ref{Sec3}.
\end{remark}

\begin{example}[A nontransitive profile on the rational branch] \label{ex:branch}
Take
\[
 v=\left(\frac9{32},\frac{27}{400},\frac{113}{800},
          \frac{43}{400},\frac{171}{800},\frac{151}{800}\right).
\]
The six components are positive and sum to one. Here
\[
 d=\frac9{16},\qquad L=\frac{279}{800},\qquad K=\frac{81}{1600},
\]
and Theorem~\ref{thm:branch} gives
\[
 p=\frac15,\quad q=\frac1{20},\quad h=\frac34,\quad
 r=0,\quad s=\frac34,\quad
 t=\frac12,\quad u=\frac3{25},\quad 1-t-u=\frac{19}{50}.
\]
These probabilities are admissible and reproduce the displayed profile by
\eqref{eq:polynomials}. Its pairwise probabilities are
\[
 \xi=\frac9{16},\qquad \eta=\frac{53}{100},\qquad
 \zeta=\frac{51}{100},
\]
all strictly larger than $1/2$. Thus the semialgebraic region contains
nontransitive profiles as well as general ranking distributions.
\end{example}

\subsection{Scope of the support restriction}
The full image of \eqref{eq:polynomials} under \eqref{eq:simplices} is the
projection of a semialgebraic set onto the $v$-coordinates. Consequently, it
is itself semialgebraic, and quantifier elimination can in principle produce
a necessary-and-sufficient description. Theorem~\ref{thm:branch} gives such
a description only for the specified branch, with the stated positivity
assumptions. We do not know a concise intrinsic characterization of the
image of the unrestricted six-dimensional parameter domain. Nor do we claim
that these particular interlaced supports realize every profile that is
realizable by independent random variables.

\section{Conclusions} \label{Sec7}

The complete six-order distribution of the classical Condorcet profile can be realized by three independent weighted dice with total support size seven, and this number is minimal. The obstruction to total support size six uses the zero ordering probability together with the positive pairwise marginals. For arbitrary ranking distributions, pairwise feasibility is governed by Trybu\l a's theorem, but it does not imply realizability of a prescribed full profile. The chain bounds provide additional necessary conditions that involve the full ordering vector. For interlaced three-point potential supports, realizability is a polynomial feasibility problem. On one boundary branch we obtain a short, exact polynomial criterion and an explicit rational construction. Natural further questions concern the other branches, boundary profiles with vanishing ordering probabilities, and bounds on the total support size needed for general realizable profiles.

\section*{Acknowledgements}
\addcontentsline{toc}{section}{Acknowledgements}

The research was supported by the National Research, Development and Innovation Office under Grants Advanced 152220 and FK 145838, and the J\'anos Bolyai Research Scholarship of the Hungarian Academy of Sciences.

\end{document}